\magnification=\magstep1       
\hsize=16.5truecm                     
\vsize=22truecm                       
\parindent 0pt                        
\parskip=\smallskipamount
\mathsurround=1pt
\voffset=2\baselineskip               
%
%
\def\today{\ifcase\month\or
  January\or February\or March\or April\or May\or June\or
  July\or August\or September\or October\or November\or December\fi
  \space\number\day, \number\year}
%
%
 at 10truept

%
%
\newcount\dispno      
\dispno=1             
\newcount\refno       
\refno=1              
\newcount\citations   
\citations=0          
\newcount\sectno      
\sectno=0             
\newbox\boxscratch    
%
%
%
%
\def\Section#1#2{\global\advance\sectno by 1\relax%
\label{Section\noexpand~\the\sectno}{#2}%
\smallskip
\goodbreak
\setbox\boxscratch=\hbox{\bf Section \the\sectno.\ }%
{\hangindent=\wd\boxscratch\hangafter=1
{\bf Section \the\sectno.\ #1}\nobreak\smallskip\nobreak}}
%
\def\sqr#1#2{{\vcenter{\vbox{\hrule height.#2pt
              \hbox{\vrule width.#2pt height#1pt \kern#1pt
              \vrule width.#2pt}
              \hrule height.#2pt}}}}
\def\square{$\mathchoice\sqr34\sqr34\sqr{2.1}3\sqr{1.5}3$}
\def\endproof{~~\hfill\square\par\medbreak}
\def\noproof{~~\hfill\square}
%
%
\def\proc#1#2#3{{\hbox{${#3 \subseteq} \kern -#1cm _{#2 /}\hskip 0.05cm $}}}

%

%
\def\normalin{\hbox{\raise0.045cm \hbox
                   {$\underline{\triangleleft }$}\hskip0.02cm}}
%
%
\def\'#1{\ifx#1i{\accent"13 \i}\else{\accent"13 #1}\fi}
%
%
%
\def\semidirect{\rlap{$\times$}\kern+7.2778pt \vrule height4.96333pt
width.5pt depth0pt\relax\;}
%
%
\def\prop#1#2{{\bf Proposition~\the\sectno.\the\dispno. }%
\label{Proposition\noexpand~\the\sectno.\the\dispno}{#1}\global\advance\dispno 
by 1{\it #2}\smallbreak}
\def\thm#1#2{{\bf Theorem~\the\sectno.\the\dispno. }%
\label{Theorem\noexpand~\the\sectno.\the\dispno}{#1}\global\advance\dispno
by 1{\it #2}\smallbreak}
\def\cor#1#2{{\bf Corollary~\the\sectno.\the\dispno. }%
\label{Corollary\noexpand~\the\sectno.\the\dispno}{#1}\global\advance\dispno by
1{\it #2}\smallbreak}
\def\defn{{\bf Definition~\the\sectno.\the\dispno. }\global\advance\dispno by 
1\relax}
\def\lemma#1#2{{\bf Lemma~\the\sectno.\the\dispno. }%
\label{Lemma\noexpand~\the\sectno.\the\dispno}{#1}\global\advance\dispno by
1{\it #2}\smallbreak}
\def\rmrk#1{{\bf Remark~\the\sectno.\the\dispno.}%
\label{Remark\noexpand~\the\sectno.\the\dispno}{#1}\global\advance\dispno
by 1\relax}

\def\proof{{\it Proof: }}

\def\numbeq#1{\the\sectno.\the\dispno\label{\the\sectno.\the\dispno}{#1}%
\global\advance\dispno by 1\relax}
\def\twist#1{\!\times_{#1}\!}
\def\comm#1,#2{\left[#1{,}#2\right]}
\newdimen\boxitsep \boxitsep=0 true pt
\newdimen\boxith \boxith=.4 true pt 
\newdimen\boxitv \boxitv=.4 true pt
\gdef\boxit#1{\vbox{\hrule height\boxith
                    \hbox{\vrule width\boxitv\kern\boxitsep
                          \vbox{\kern\boxitsep#1\kern\boxitsep}%
                          \kern\boxitsep\vrule width\boxitv}
                    \hrule height\boxith}}
\def\square{\ \hbox{\vrule height7.5pt depth1.5pt width 6pt}\par}
\outer\def\square{\ifmmode\else\hfill\fi
   \setbox0=\hbox{} \wd0=6pt \ht0=7.5pt \dp0=1.5pt
   \raise-1.5pt\hbox{\boxit{\box0}\par}
}

\def\frac#1/#2{\leavevmode\kern.1em
              \raise.5ex\hbox{\the\scriptfont0 #1}\kern-.1em
              /\kern\.15em\lower.25ex\hbox{\the\scriptfont0 #2}}
\def\incnoteq{\lower.1ex \hbox{\rlap{\raise 1ex
     \hbox{$\scriptscriptstyle\subset$}}{$\scriptscriptstyle\not=$}}}
%
%
\def\mapright#1{\smash{
     \mathop{\longrightarrow}\limits^{#1}}}
\def\mapdown#1{\Big\downarrow
 \rlap{$\vcenter{\hbox{$\scriptstyle#1$}}$}}


\def\propcontup{\bigcup\!\!\!\rlap{\kern+.2pt$\backslash$}\,\kern+1pt\vert}
%
%
%
\def\label#1#2{\immediate\write\aux%
{\noexpand\def\expandafter\noexpand\csname#2\endcsname{#1}}}
%
\def\ifundefined#1{\expandafter\ifx\csname#1\endcsname\relax}
%
%
\def\ref#1{%
\ifundefined{#1}\message{! No ref. to #1;}%
 \else\csname #1\endcsname\fi}
%
%
\def\refer#1{%
\the\refno\label{\the\refno}{#1}%
\global\advance\refno by 1\relax}
%
%
\def\cite#1{%
\expandafter\gdef\csname x#1\endcsname{1}%
\global\advance\citations by 1\relax
\ifundefined{#1}\message{! No ref. to #1;}%
\else\csname #1\endcsname\fi}
%
%
\font\bb=msbm10 
 at 8truept      
%
%
%

\def\Q{\hbox{\bb Q}}

\def\Z{\hbox{\bb Z}}

\def\Z{\hbox{\bb Z}}                     

\newread\aux
\immediate\openin\aux=\jobname.aux
\ifeof\aux \message{! No file \jobname.aux;}
\else \input \jobname.aux \immediate\closein\aux \fi
\newwrite\aux
\immediate\openout\aux=\jobname.aux

\font\smallheadfont=cmr8 at 8truept

\headline={\ifnum\pageno<2{\hfill}\else{\smallheadfont
Bilinear Maps and Central Extensions \hfill Arturo
Magidin}\fi}

\centerline{\bf BILINEAR MAPS AND CENTRAL EXTENSIONS OF ABELIAN GROUPS}
\centerline{Arturo Magidin\footnote*{The author was 
supported by a fellowship from the Programa de Formaci\'on y
Superaci\'on del Personal Acad\'emico de la UNAM, administered by the
DGAPA.}}

\centerline{Last modified: \today}
\smallskip
{\parindent=20pt
\baselineskip=12pt
\narrower\narrower
\noindent{\smallheadfont {\smallheadfont Abstract.\/} 
Every nilpotent group of class at most two may be 
embedded in a central extension of abelian groups with bilinear
cocycle. The embedding is shown to depend only on the base 
group. Some refinements are obtained by considering the
cohomological situation explicitly.}\par\smallskip}

\footnote{}{\smallheadfont AMS Classification:20E22, 20F18 (primary);
20J05 (secondary). Keywords:central extensions,nilpotent of class two.}
\bigbreak

The main result of this paper is that every nilpotent group of class
at most two may be embedded into a central extension of abelian
groups, in which the associated cocycle is~bilinear (definitions are
recalled in \ref{prelims}). The result is
related to a paper of N.J.S.~Hughes (see~{\bf [\cite{hughes}]}), in
which he establishes a one to one correspondence between the
equivalence classes of central extensions of abelian groups with what
he calls {\it ``regular bilinear mappings''}. The proof of Theorem~2
in {\bf [\cite{hughes}]} is very similar to our proof of
\ref{fullextension} in the present~work.

In the first section we will recall the basic definitions relating to
central extensions and the second cohomology group of two abelian
groups. We will also introduce the notion of a twisted product of
abelian groups, to be used later. In the second section we establish
our main result and discuss some examples. In the third and final
section we cast the main theorem into a cohomological setting to gain
more information about the embedding we construct in \ref{centextandtwists}

It is my very great pleasure to express my deep gratitude and
indebtedness to Prof.~George M.~Bergman, for his ever ready advice and
encouragement throughout the preparation of this~paper.
\smallskip

\Section{Preliminaries and notation}{prelims}

We will denote the identity element of a group~$G$ by $e_G$ if~$G$ is
written multiplicatively, and by $0_G$ if it is written additively; we
will omit the subscript if it is clear from context. Given a group~$G$
and elements $x,y\in G$, $[x,y]$ denotes their commutator; that is
$$[x,y]=x^{-1}y^{-1}xy.$$ Given two subgroups $U$ and~$V$ of $G$,
$[U,V]$ denotes the subgroup of~$G$ generated by all elements $[u,v]$,
where $u\in U$ and $~v\in V$.

For general notation and a more complete treatment of the extension
problem, we direct the reader to {\bf [\cite{rotman}]}, Chapter~7. We
will also omit proofs and direct the reader to Rotman's~text for them.  Let
$$1\,\mapright{}\,B\,\mapright{i}\,G\,\mapright{\pi}\,
A\mapright{}\,1$$ be a central extension (so that $B\subseteq Z(G)$)
of groups with $B$ and~$A$, both abelian.  It will be convenient (and
the formulas clearer) if the product in~$G$ (and hence in~$B$) is
written additively, and the product in~$A$ multiplicatively, although
in general $G$ may be~nonabelian.

Recall that a {\it factor set} (or {\it cocycle}) is a map
$\gamma\colon A\times A\to B$ that satisfies the identities:
$$\displaylines{
\qquad(\numbeq{cocycleone})\;\forall\,x,y\in
A,\>\gamma(e,y)=\gamma(x,e) = 0\hfill\cr
\qquad(\numbeq{cocycletwo})\;\hbox{The {\it cocycle identity} holds for
every $x$, $y$, and $z$ in $A$:}\hfill\cr
\gamma(x,y) + \gamma(xy,z) = \gamma(y,z)+\gamma(x,yz).\cr}$$

Any central extension determines a factor set $\gamma$, once a transversal
$\ell\colon A\to G$ with $\ell(e)=0$ and $\pi\circ\ell = {\rm id}_A$
has been chosen, by letting
$\gamma(x,y)$ be the element of~$B$ such that $$\ell(x)+\ell(y)=
\ell(xy) + \gamma(x,y).$$  
Note that $\ell$ need not be (and in general will not be) a group
morphism.

Conversely any factor set $\delta$ satisfying (\ref{cocycleone})
and~(\ref{cocycletwo}) determines a central extension and a
transversal. Namely, we take the set of all pairs $(a,b)$ with $a\in
A$ and~$b\in B$ as the underlying set for the extension, and define
the multiplication by the rule
$$(a,b)\cdot (a',b') = \bigl(aa', b+b'+\delta(a,a')\bigr)$$
and the transversal by $\ell(a)=(a,0)$.

We note that if $\delta$ is a bilinear map from $A\times A$ to~$B$, then it
satisfies (\ref{cocycleone}) and~(\ref{cocycletwo}), although
the converse is not necessarily true. If the factor
set $\delta$ is bilinear, we will denote the resulting group by
$A\twist\delta B$ and say that this group is a {\it twisted product
of~$B$ by~$A$}.

A function $\eta\colon A\times A \to B$ is called a {\it coboundary}
iff there is a set map $h\colon A \to B$ with $h(e)=0$ and such that
for all $x$ and $y$ in~$A$,
$$\eta(x,y)=h(x) + h(y) - h(xy).$$

We denote the set of all factor sets by $Z^2(A,B)$, and the set of all
coboundaries by $B^2(A,B)$. They both carry a structure of abelian
groups, given by pointwise addition, and $B^2(A,B)$ is a subgroup of
$Z^2(A,B)$. We define the {\it second cohomology group of $A$ with
coefficients in~$B$} to be the quotient group
$$H^2(A,B) \cong Z^2(A,B)/B^2(A,B).$$

It is not hard to verify that given a map $\phi\colon B\to C$, we
obtain an induced map $$\phi^*\colon H^2(A,B)\to H^2(A,C)$$
by taking the class of the factor set $\delta\colon A\times A\to B$ to
$\phi\circ\delta$. Dually, a map $\psi\colon A\to C$ induces a map
$$\psi_*\colon H^2(C,B)\to H^2(A,B)$$
by taking the class of the factor
set $\gamma\colon C\times C\to B$ to the factor set $\gamma\circ
(\psi\times\psi)$.

Recall also that two central extentions of~$B$ by~$A$,
$1\>\mapright{}\>B\>\mapright{i}\>G\>\mapright{\pi}\>A\>
\mapright\>1$~and
$1\>\mapright{}\>B\>\mapright{i'}\>G'\>\mapright{\pi'}\>
A\>\mapright{}\>1$ are said to be {\it equivalent}
iff there exists a homomorphism
\hbox{$\phi\colon G\to G'$} making the following diagram of exact sequences
commute:
{$$
\matrix{1&\mapright{}&B&\mapright{i}&G&
\mapright{\pi}&A&\mapright{}&1\cr
&&\|&&\;\mapdown{\phi}&&\|\cr
1&\mapright{}&B&\mapright{i'}&G'&\mapright{\pi'}
&A&\mapright{}&1\cr}$$}

By a classical theorem of Schreier, the elements of $H^2(A,B)$ are in
one to one correspondence with the equivalence classes of central
extensions of $B$ by~$A$. The extensions which are abelian groups form
a subgroup, namely ${\rm Ext}(A,B)$.

We will denote the set of all equivalence classes of central
extensions of~$B$ by~$A$, at least one of whose representatives is a
twisted product, by $H^2_{\rm Bil}(A,B)$. Note that
$$\eqalign{H^2_{\rm Bil}(A,B)&\cong \Bigl({\rm Hom}(A\otimes A,B) +
B^2(A,B)\Bigr)\Bigm/B^2(A,B)\cr
&\cong {\rm Hom}(A\otimes A,B)\Bigm/ \Bigl({\rm Hom}(A\otimes A,B)\cap
B^2(A,B)\Bigr)\cr}$$

\rmrk{nottwist} There are central extensions which are not
equivalent to a twisted product; that is,
extensions for which the
factor set associated to a transversal does not lie in $H^2_{\rm
Bil}(A,B)$. For example, let $p$ be an odd prime, and consider the
extension $$0\,\mapright{}\, \Z/p\Z\,
\mapright{i}\,\Z/p^2\Z\,\mapright{\pi}
\,\Z/p\Z\,\mapright{}\,0$$
given by $i(\overline{1})=\overline{p}$ (we are viewing
$\Z/a\Z$ as the integers modulo $a$, so it has the class of $1$ as a
distinguished generator; also note that we are writing all groups
additively in this~example).

As a transversal $\ell\colon \Z/p\Z\to \Z/p^2\Z$, we can take the map
$\ell(\overline{x})=\overline{x}$ for $0\leq x<p$. Then, the associated factor
set $\gamma$ must satisfy
$$\ell(\overline{x}) + \ell(\overline{y}) = \ell(\overline{x+y}) +
\gamma(\overline{x},\overline{y})$$ for every $0\leq x,y <
p$. Therefore, we have
$$\gamma(\overline{x},\overline{y}) = \cases{0&if $0\leq x+y < p$;\cr
p&if $p\leq x+y$\cr}$$
for all $x$ and~$y$, with $0\leq x,y < p$. This is clearly not
bilinear, since, for example $\gamma(\overline{p-1},\overline{p-1}) =
\overline{p}$, but $\gamma(0,0) - 
\gamma(0,\overline{1}) - \gamma(\overline{1},0) +
\gamma(\overline{1},\overline{1}) = 0$.

Is $\gamma$ cohomologous to a bilinear map? The answer is no. If it
were, we would be able to write $\Z/p^2\Z$ as a twisted product
$\Z/p\Z\twist\delta\Z/p\Z$; an easy calculation will show that for
every $(x,y)$ in such a group,
$$\eqalign{(x,y)^p &= \left(px, py + {p(p-1)\over 2}\delta(x,x)\right)\cr
&= (0,0)\cr}$$
contradicting our assumption that such a twisted product is isomorphic to
$\Z/p^2\Z$. So $\gamma$ is not cohomologous to a bilinear map.

\bigbreak
\Section{Central extensions and twisted products}{centextandtwists}

We are now ready to describe our result. Let $A$ and $B$ be two
abelian groups, and let~$G$ be a central extension of~$B$ by~$A$, with
transversal~$\ell$ and factor set~$\gamma$. We will show that there exists an
abelian group~$L$, an embedding $\varphi\colon B\to L$, a twisted
product~$G'$ of~$L$ by~$A$, and an embedding $\psi\colon G\to G'$ such
that the following diagram of exact rows commutes:

{$$
\matrix{1&\mapright{}&B&\mapright{}&G&\mapright{}&A&\mapright{}&1\cr
&&\mapdown{\varphi}&&\mapdown{\psi}&&\|\cr
1&\mapright{}&L&\mapright{}&G'&\mapright{}&A&\mapright{}&1\cr}$$}

This will imply that the twisted product $G'$ realizes the
factor set $\varphi^*(\gamma)$. That is, that $\varphi^*(\gamma)$ is
cohomologous to the bilinear factor associated to the
twisted product~$G'$.

In fact, the group $L$ may be taken to depend {\sl only on~$B$}, and
not on the group~$A$ or the factor set~$\gamma$. In other
words, given an abelian group~$B$ there exists an abelian group~$L$
and an embedding $\varphi\colon B\to L$ such that for all abelian
groups~$A$, the induced map
$$\varphi^*\colon H^2(A,B)\, \mapright{}\, H^2(A,L)$$
has image in~$H^2_{\rm Bil}(A,L)$.

Recall that a bilinear map $\alpha\colon A\times A\to B$ is said to be
{\it alternating} if $\alpha(x,x)=0$ for all $x\in A$. It is not hard
to verify that this condition implies that $\alpha(x,y)=-\alpha(y,x)$,
and that the converse holds if $B$ has no $2$-torsion.

Suppose that we have two abelian groups~$A$ and~$B$, and an
alternating bilinear map $\alpha\colon A\times A\to B$. Let the
abelian group~$C$, the bilinear map $\beta\colon A\times A\to C$, and
the map $i_B\colon B\to C$ be the universal triple with the following
universal property:

{\parindent=10pt
 \narrower\narrower
\item{(\numbeq{conditionone})}{The map $i_B\colon B\to C$ makes
 the following diagram commute:
$$\matrix{&C\hfill\cr
          \hfill{{\scriptscriptstyle\beta(x,y)-\beta(y,x)}\atop{}}%
\!\!\!\!\nearrow&\cr
          A\times A&\hfill{\Big\uparrow}i_B\cr
          \hfill{{}\atop{\scriptscriptstyle
          \alpha(x,y)}}\!\!\!\!\!\searrow\cr
          &B\hfill\cr}$$}

\item{(\numbeq{conditiontwo})}{For any abelian group $C'$, any bilinear map
$\beta'\colon
A\times A \to C'$, and map $\varphi'\colon B\to C'$ such that 
$$\matrix{&C'\hfill\cr
          \hfill{{\scriptscriptstyle\beta'(x,y)-\beta'(y,x)}\atop{}}%
\!\!\!\!\nearrow&\cr
          A\times A&\hfill{\Big\uparrow}\varphi'\cr
          \hfill{{}\atop{\scriptscriptstyle
          \alpha(x,y)}}\!\!\!\!\!\searrow\cr
          &B\hfill\cr}$$
commutes, there exists a unique $\psi\colon C\to C'$ such that
$\varphi'=\psi\circ i_B$ and $\beta'=\psi\circ\beta$.}}

\rmrk{injectivity} We claim that the universal map $i_B\colon B\to
C$ must be injective.
For this, first we examine the case of $A$ finitely generated, say
$$A\cong\Z^r\oplus \Z/a_{r+1}\Z\oplus\cdots\oplus\Z/a_{r+m}\Z$$
with $a_{r+1}|\cdots|a_{r+m}$ (we write $A$ additively for the
remainder of this remark). Let $g_1,\ldots,g_{r+m}$ be the
canonical set of generators for $A$ under this isomorphism, and take 
$x,y\in A$, given by
$$x = k_1g_1 + \cdots +k_{r+m}g_{r+m}\leqno(\numbeq{expx})$$
$$y = l_1g_1 + \cdots +l_{r+m}g_{r+m}\leqno(\numbeq{expy})$$
where $k_1,\ldots,k_{r+m},l_1,\ldots,l_{r+m} \in \Z$ and $0\leq
k_{r+i},l_{r+i}< a_{r+i}$ for $i=1,\ldots,m$. We let $C'=B$ and
define
$$\beta'(x,y) = \prod_{1\leq i<j\leq
r+m}{k_jl_i}\Bigl(\alpha(g_i,g_j)\Bigr)$$

It is now easy to see that for $x,y\in A$ as in (\ref{expx})
and (\ref{expy}), we have $$\beta'(x,y)-\beta'(y,x) = \alpha(x,y)$$ so
we can set up a commutative diagram as in (\ref{conditiontwo}) with
$B$ playing the role of~$C'$; by the universal property of the triple
$(C,i_B,\beta)$, there exists a unique map $\psi\colon C \to
B$ with ${\rm id}_B =\psi\circ i_B$. So, in particular, $i_B$
must be injective. 

For the general case, we consider finitely generated subgroups $A'$
of~$A$, and the restriction of $\alpha$ to $A'\times A'$. The triples
$(A',B,\alpha|_{A'\times A'})$ form a
directed system, with the connecting maps being the
obvious inclusions. The universal object $C$ will be the limit of the
universal objects $C'$ at each level, and the universal map $i_B$
will be the direct limit of the maps at each level. Since direct
limits respect monomorphisms, the universal map $i_B\colon B\to C$
is injective, establishing the~claim. 

Let $G$ be a central extension of $B$ by $A$, with exact sequence 
$$1\,\mapright{}\, B\mapright{i}\,G\,\mapright{\pi}\, A\,\mapright{}\,
1,
\leqno(\numbeq{exseq})$$
let $L$ an abelian group (written additively), and $\beta\colon
A\times A\to L$ a bilinear map. We abuse notation slightly
and denote the induced bilinear map on $G\times G$ also by $\beta$. We
denote the image of an arbitrary $g\in G$ under $\pi$ by $\overline{g}$.

We want to know what conditions a set map $f\colon G\to L$ must
satisfy in order for the map $\varphi\colon G\to A\twist\beta L$,
given by $$\varphi(g)=\Bigl(\overline{g},f(g)\Bigr),\leqno(\numbeq{vrphi})$$ to
be an injective group morphism.

\lemma{conditionsformap}{Let $G$,~$L$ and~$\beta$ be as in
(\ref{exseq}) and the subsequent two paragraphs. A set map
$f\colon G\to L$ makes $\varphi$, defined as in (\ref{vrphi}), a
group morphism if and only if
$$\forall x,y\in G \qquad f(x+y) = f(x) + f(y) +
\beta(x,y).\leqno(\numbeq{mapcondone})$$
Furthermore, such an $f$ makes $\varphi$ into an injective group
morphism if and only if it also satisfies
$$f|_{B}\colon B \to L\qquad \hbox{is a one-to-one group
morphism}.\leqno(\numbeq{mapcondtwo})$$}

\proof Let $x,y\in G$. Then we must have
$$\eqalign{(\overline{x},f(x))\cdot(\overline{y},f(y)) &=
       (\overline{x}\,\overline{y},f(x) + f(y) + \beta(x,y))\cr
&=(\overline{x+y}, f(x+y))\cr}$$
so condition (\ref{mapcondone}) is both necessary and sufficient.

Assuming we already have condition (\ref{mapcondone}), we want to know
when
$\varphi$ is injective. Clearly, if $x$ and $y$ are not congruent
modulo $B$, then $\overline{x}\neq\overline{y}$ in $A$, so
$\varphi(x)\neq\varphi(y)$ is guaranteed. So it is both necessary
and sufficient to require that $f$ be an injective set map on
$B$. Since $B$ lies in both the kernel on the right
and on the left of $\beta$, this condition is equivalent to asking
that $f|_B$ be an injective group morphism, giving
condition (\ref{mapcondtwo}).\endproof

\rmrk{powerform} If $f$ satisfies condition (\ref{mapcondone}), then
$$\forall x\in G,\>\forall n\in\Z\quad f(nx) = nf(x) +
{n\choose2}\beta(x,x)$$ where ${n\choose2}= {n(n-1)\over 2}$ for all
$n\in\Z$. This is easily established by induction on~$n$.

\rmrk{commrmrk} An easy calculation
shows that if $f$ satisfies~(\ref{mapcondone}) then $$\forall x,y\in
G\qquad\quad f(\comm x,y) =
\beta(x,y) - \beta(y,x).\leqno(\numbeq{commformula})$$

We want to show that for $G$ as in (\ref{exseq}) and the subsequent
two paragraphs, we can always find a group~$L$ and define a set map
$f$ satisfying conditions (\ref{mapcondone}) and
(\ref{mapcondtwo}). We will do this by defining $f$ on the subgroup
$B$, making sure it satisfies (\ref{commformula}), and then showing it
can be extended to all of $G$. For this, we need an
extension result.

\lemma{extension}{Let $G$ be as above, and let $U$ and~$V$ be two
subgroups of~$G$,
with \hbox{$[V,U]\subseteq U$}. Let $L$ be an abelian group, and
$$\beta\colon A\times A\to L$$ be a bilinear map. Suppose
we are given functions
$$\eqalign{f_U&\colon U\to L\cr
           f_V&\colon V\to L\cr}$$
with $f_U|_{\scriptscriptstyle U\cap V} = f_V|_{\scriptscriptstyle
V\cap U}$, and satisfying (\ref{mapcondone}) (That is, $f_U$ satisfies
(\ref{mapcondone}) with $x$ and~$y$
ranging over $U$, and $f_V$ with $x$ and $y$ ranging over
$V$). Suppose further that $f_U$ satisfies equation 
(\ref{commformula}) whenever
$x\in V$ and $y\in U$. Let $W=\langle U,V\rangle=UV$ be the subgroup
generated by $U$ and $V$. Finally, let $f_W\colon W\to L$ be given by
$$f_W(u+v) = f_U(u) + f_V(v) + \beta(u,v).$$
Then $f_W$ is an extension of $f_U$ and $f_V$, and satisfies
(\ref{mapcondone}) (with
$x$ and~$y$ ranging over~$W$).}

\proof  Note that since $[G,G]\subseteq B$, and $B$ lies in both the
kernel on the right and the kernel on the left of~$\beta$, it follows
that $[U,V]$ also lies in both the kernel on the right
and kernel on the left of~$\beta$.

We need to check that $f_W$ is well defined, that it indeed
extends both $f_U$ and $f_V$, and that it satisfies
(\ref{mapcondone}).

Suppose that $u_1+v_1 = u_2+v_2$. Then $-u_2+u_1 = v_2-v_1\;\in
U\cap V$. So in particular we must have that $f_U(-u_2+u_1) =
f_V(v_2-v_1)$, i{.}e{.}
$$f_U(u_1) - f_U(u_2) + \beta(u_2,u_2) -\beta(u_2,u_1) =
  f_V(v_2) - f_V(v_1) + \beta(v_1,v_1) - \beta (v_2,v_1).$$

We also have $u_1=u_2+v_2-v_1$ and $v_2 = -u_2+u_1+v_1$. So,
according to the definitions, we have
$$\eqalign{\scriptstyle f_W(u_1+v_1) - f_W(u_2+v_2) &\scriptstyle=
          f_U(u_1) + f_V(v_1) + \beta(u_1,v_1) -
                   \bigl(f_U(u_2) + f_V(v_2) + \beta(u_2,v_2)\bigr)\cr
   &\scriptstyle= f_U(u_1) - f_U(u_2) - \bigl(f_V(v_2) - f_V(v_1)\bigr) +
      \beta(u_1,v_1) - \beta(u_2,v_2)\cr
   & \scriptstyle= \beta(v_1,v_1) - \beta(v_2,v_1) -\beta(u_2,u_2) +
          \beta(u_2,u_1)+
\beta(u_2+v_2-v_1,v_1) - \beta(u_2,-u_2+u_1+v_1)\cr
   &\scriptstyle= \beta(v_1,v_1) - \beta(v_2,v_1) - \beta(u_2,u_2)+
          \beta(u_2,u_1)+ \beta(u_2,v_1) +\cr
&\scriptstyle\qquad + \beta(v_2,v_1) - \beta(v_1,v_1)
      +\beta(u_2,u_2) - \beta(u_2,u_1) - \beta(u_2,v_1)\cr
   &\scriptstyle=\  0\cr}$$
so $f_W$ is well defined.

That $f_W$ extends both $f_U$ and $f_V$ now follows
trivially. Finally,
let $u_1+v_1$ and~$u_2+v_2$ be two elements of $W$.  Then, since
$[U,V]\subseteq Z(G)$,
$$(u_1+v_1)+(u_2+v_2) = (u_1+u_2+[v_1,u_2])+(v_1+v_2)$$
with $u_1+u_2+[v_1,u_2]\in U$ and $v_1+v_2\in V$. So
$$\eqalign{\scriptstyle f_W((u_1+v_1)+ (u_2+v_2))
&\scriptstyle= f_W(u_1+u_2+[v_1,u_2]+v_1+v_2)\cr
     &\scriptstyle= f_U(u_1+u_2+[v_1,u_2]) + f_V(v_1+v_2) +
        \beta(u_1+u_2+[v_1,u_2],v_1+v_2)\cr
     &\scriptstyle= f_U(u_1) + f_U(u_2) +
     \beta(v_1,u_2)-\beta(u_2,v_1)+\beta(u_1,u_2)+
      f_V(v_1) + f_V(v_2) + \beta(v_1,v_2)+\cr
     &\scriptstyle\quad +\beta(u_1,v_1)
     +\beta(u_1,v_2)+\beta(u_2,v_1)+\beta(u_2,v_2)\cr
     &\scriptstyle= f_U(u_1) + f_U(u_2) + \beta(v_1,u_2) + \beta(u_1,u_2) +
     f_V(v_1) + f_V(v_2)+\beta(v_1,v_2)+\beta(u_1,v_1)+\cr
     &\scriptstyle\quad + \beta(u_1,v_2) + \beta(u_2,v_2)\cr}$$

On the other hand we have
$$\eqalign{\scriptstyle f_W(u_1+v_1) + f_W(u_2+v_2) +
\beta(u_1+v_1,u_2+v_2) 
&\scriptstyle= f_U(u_1) + f_V(v_1) + \beta(u_1,v_1) +
   f_U(u_2) + f_V(v_2) + \cr 
   &\scriptstyle\quad +\beta(u_2,v_2) + \beta(u_1,u_2) + \beta(u_1,v_2) +
   \beta(v_1,u_2) + \beta(v_1,v_2)\cr}$$
so $f_W$ satisfies (\ref{mapcondone}) and we are done.\endproof

Now we set up some notation. Let $G$ be a central extension of~$B$
by~$A$, and $H$ a proper
subgroup of $G$ with $B\subseteq H$. Let $L$ be an abelian
divisible group, $\beta\colon A\times A\to L$ a bilinear
map, and $f\colon H\to L$ be a set map satisfying
(\ref{mapcondone}),~(\ref{mapcondtwo}), and~(\ref{commformula}).  Note
that condition (\ref{mapcondone}) is assumed only for $x$ and~$y$ in $H$,
whereas condition (\ref{commformula}) is assumed for any commutator
$[x{,}y]$ with $x$ and~$y$ arbitrary elements of $G$. Finally, let
$g\in G\backslash H$ and let
$K=\langle H,g\rangle$.

We want to extend $f$ to all of $K$, retaining properties
(\ref{mapcondone}) (this time for all $x$ and~$y$ in~$K$),
(\ref{mapcondtwo}), and~(\ref{commformula}).

\thm{fullextension}{Let $G$, $H$, $\beta$, $L$, $f$ and~$K$ be as in
the previous two paragraph. Then $f$ can be extended to $K$.}

\proof We want to apply \ref{extension}, so we let $H$ play the
role of $U$ and the cyclic group generated by $g$ play the role of
$V$. $f_U$ will simply be our given $f$, and we need to define~$f_V$.
Let $I=\{a\in\Z_+ \,|\, ag\in U\}$. We have two cases.

{\it Case 1.} $I=\emptyset$. Note in particular that in this case, $g$
cannot be a torsion element of~$G$. We define $f_V(g)=0$ and, extend
to every power of~$g$ using the formula in \ref{powerform}. Then
\ref{extension} gives us an extension of~$f$ to all of~$K$.

{\it Case 2.} $I\not=\emptyset$. Let $n_0=min\{n\in I\}$. Then $U\cap
V=\langle n_0g\rangle$. Using \ref{powerform}, we know that whatever
we define $f_V(g)$ to be, it must satisfy
$$f_V(n_0g) = n_0f_V(g) + {n_0\choose 2}\beta(g,g).$$
However, $n_0g\in U$, so $f_V(n_0g)=f_U(n_0g)$. Solving for $f_V(g)$,
we see that $f_V(g)$ must satisfy
$$n_0f_V(g) = f_U(n_0g) - {n_0\choose 2}\beta(g,g);$$
so we define
$$f_V(g) = {1\over n_0}\Bigl(f_U(n_0g) - {n_0\choose 2}\beta(g,g)\Bigr)$$
where by ${1\over n}(\cdots)$ we mean any $n$-th root, chosen now once
and for all. Then we extend $f_V$ to all of~$V$ using the formula in
\ref{powerform}. We need to verify that this is well defined if~$g$ is
a torsion element, and that it agrees with~$f_U$ on $U\cap V$. 

First, suppose that~$g$ is a torsion element. By definition of~$n_0$,
it follows that the order of~$g$ is a multiple of~$n_0$, say
$k_0n_0$. We want to verify that for every~$a\in\Z$, \hbox{$f_V(ag)=
f_V((k_0n_0+a)g)$}. 

By definition,
$$\eqalignno{f_V(ag) &= af_V(g) + {a\choose 2}\beta(g,g)\cr
&= a\left({1\over n_0}\Bigl(f_U(n_0g) - {n_0\choose
2}\beta(g,g)\Bigr)\right) + {a\choose 2}\beta(g,g).\cr}$$

On the other hand, we have
$$\eqalignno{f_V\bigl((k_0n_0+a)g\bigr) &= (k_0n_0+a)f_V(g) + {k_0n_0+a\choose
2}\beta(g,g)\cr
&=(k_0n_0+a)\left({1\over n_0}\Bigl(f_U(n_0g)-{n_0\choose
2}\beta(g,g)\Bigr)\right) + {k_0n_0+a\choose 2}\beta(g,g)\cr}$$

A long but straightforward calculation now shows that
the values indeed agree.

As for the values on the intersection of~$U$ and~$V$, since the
intersection is generated by $n_0g$, we need to check the values
of~$f_V$ at $kn_0g$, for $k\in\Z$. We~have
$$\eqalignno{f_V\bigl((kn_0)g\bigr)&= kn_0f_V(g) + {kn_0\choose 2}\beta(g,g)\cr
     &= k\left(f_U(n_0g) - {n_0\choose 2}\beta(g,g)\right)
        + {kn_0\choose 2}\beta(g,g)\cr
     &= kf_U(n_0g) + \left( {kn_0\choose 2} - k{n_0\choose
       2}\right)\beta(g,g)\cr
     &= kf_U(n_0g) + \left( k{n_0\choose 2} + {k\choose 2}n_0^2 -
       k{n_0\choose 2}\right)\beta(g,g)\cr
     &= kf_U(n_0g) + {k\choose 2}n_0^2\beta(g,g)\cr
     &= kf_U(n_0g) + {k\choose 2}\beta(n_0g,n_0g)\cr
     &= f_U\left(k(n_0g)\right) \cr
     &= f_U((n_0k)g)\cr}$$
so $f_U|_{U\cap V} = f_V|_{V\cap U}$. Applying \ref{extension}, we get
an extension of~$f$ to all of~$K$.\endproof

\rmrk{divisibility} Notice that the divisibility of $L$ was only used
in defining $f_V(g)$, in the case where $I$ is nonempty. So, for
example, if we know that the quotient $G/H$ has exponent $k$, then we
would only need $L$ to have all $k$-th roots.

\thm{centextthentwist}{Let $1\,\mapright{}\, B\,\mapright{i}\, G\,
\mapright{\pi}\, A\,\mapright\,{}1$ be a central extension of abelian
groups. Then there exists an abelian group~$L$, a bilinear map
\hbox{$\beta\colon A\times A\to L$} and a set map \hbox{$f\colon G\to
L$} satisfying (\ref{mapcondone}) and (\ref{mapcondtwo}).}

\proof  For the moment assume $L$ given. We will consider the
set of pairs $(H,f_H)$, where $H$ is a subgroup of $G$ containing~$B$,
and $f_H\colon H\to L$ is a set function satisfying
(\ref{mapcondone}), (\ref{mapcondtwo}), and~(\ref{commformula}). We
partially order the set by inclusion on the first coordinate, and
``extension'' on the second. i{.}e{.} $(H,f_H)\leq(H',f_H')$ iff
$H\subseteq H'$ and $f_{H'}|_{\scriptscriptstyle H} = f_H$.

We need to verify that there is an $L$ for which this set of pairs is
non-empty. For that, it suffices to show that there exists an abelian
group $L$, a bilinear map $\beta\colon A\times A\to L$,
and a set map $f\colon B\to L$ satisfying (\ref{mapcondtwo})
and~(\ref{commformula}) (for $x,y\in B$, (\ref{mapcondone}) follows
from~(\ref{mapcondtwo})).

To establish the existence of~$L$, we consider the alternating
bilinear map $\alpha(x,y)=[x,y]$.  Then \ref{injectivity}
shows that there is an abelian group~$C$, an embedding
\hbox{$i_B\colon B\to C$}, and a bilinear map
$\beta\colon A\times A\to C$ satisfying $\beta(x,y)-\beta(y,x) =
[x,y]$ for all~$x$ and~$y$ in~$G$, with the triple $(C,i_B,\beta)$
being~universal.

Let $L$ be a divisible abelian group containing~$B$, and $j\colon B\to
L$ an embedding. Since~$L$ is divisible, it is an injective object of
the category of abelian groups, so the inclusion $i_B\colon B\to C$
factors through $j$; that is, there exists a map $\chi\colon C\to L$
(not necessarily a monomorphism), such that $\chi\circ i_B=j$. We now
let $\beta'=\chi\circ\beta$, and we let $f\colon B\to L$ be equal to
$j$. Clearly, $f$ satisfies~(\ref{mapcondtwo}). Also, 
$$\eqalign{f([x,y])&=j([x,y])\cr
&= j(\alpha(x,y))\cr
&= \chi\bigl( i_B(\alpha(x,y))\bigr)\cr
&= \chi\bigl(\beta(x,y) - \beta(y,x)\bigr)\cr
&= \beta'(x,y) - \beta'(y,x),\cr}$$
so $f$ satisfies (\ref{mapcondtwo})
and~(\ref{commformula}). Therefore, the set of pairs defined in this
first paragraph of this proof is non-empty, for the given~$L$.

Finally, we apply Zorn's Lemma to this set of
pairs. \ref{fullextension} implies that a maximal element in our set
must have first coordinate equal to~$G$, so we are done.\endproof

\cor{embedcentext}{Let $1\,\mapright{}\, B\,\mapright{i}\, G\,
\mapright{\pi}\, A\,\mapright{}\, 1$ be a central extension of $B$ by
$A$. Then there exists an abelian group~$L$ containing~$B$, a bilinear
map $\beta\colon A\times A \to L$,
and an embedding $\varphi\colon G\to A\twist\beta L$ such that the
following diagram commutes and has exact rows:
{$$
\matrix{1&\mapright{}&B&\mapright{i}&G&\mapright{\pi}&A&\mapright{}&1\cr
&&\mapdown{j}&&\mapdown{\varphi}&&\|\cr
1&\mapright{}&L&\mapright{i'}&A\twist\beta
L&\mapright{\pi'}&A&\mapright{}&1\cr}$$}
Furthermore, $L$ depends only on~$B$.}

\proof Let $L$, $f$ and $\beta$ be as in \ref{centextthentwist}. We
note that the group~$L$ we obtain this way depends only on~$B$, and
not on the groups $A$ or~$G$. The
map $\varphi$ is simply given by the projection onto $A$ in
the first coordinate, and the map~$f$ in
the second coordinate. It only remains to show that the diagram
commutes.

Because the first component of $\varphi$ is just $\pi$, commutativity
of the square on the right is immediate. For the square on the left,
let $b\in B$. Then $\varphi(i(b))= (1,b)$ by construction of~$f$,
which is also the same as $j(i'(b))$, since $j$ is an embedding and
$i'$ is just the inclusion of~$L$ into $A\twist\beta L$. So the diagram
commutes and we are done.\endproof

{\bf Example \numbeq{niltwoproduct}.} Let $p$ be an odd prime, and
let~$G$, written multiplicatively, be the group presented by
$$G = \Bigl\langle x,y,z\,\Bigm|\, x^p=y^p=[x,y]^p=[[x,y],x] =
[[x,y],y] = [x,z]=[y,z] = e;\; z^p=[x,y]\Bigr\rangle.$$
Let~$B$ be the subgroup of~$G$ generated by $[x,y]$. Then $B=[G,G]$ and
$A\cong G/B$ is the abelianization of~$G$, that is
$$A= \Bigl\langle \overline{x},\overline{y},\overline{z}\,\Bigm|\,
\overline{x}^p=\overline{y}^p=\overline{z}^p=
[\overline{x},\overline{y}]= [\overline{x},\overline{z}] =
[\overline{y},\overline{z}]=e\Bigr\rangle$$
so~$A$ is the product of three cyclic groups of order~$p$.

Let $\beta\colon G\times G\to B$ be given by commutator map; that is,
\hbox{$\beta(g,g')=[g,g']$}. Since~$G$ is nilpotent of class two, it
is not hard to verify that~$\beta$ induces a well defined map on~$A$,
and is~bilinear. 

Is~$G$ a twisted product? The answer is no. For supposing it were,
there would exist a bilinear map \hbox{$\gamma\colon A\times A\to B$}
such that $G\cong A\times_{\gamma}B$. Since both~$A$ and~$B$ are of
exponent~$p$, and since~$p$ is odd, it follows that~$G$ is also of
exponent~$p$; but $z$ is a $p$-th power in~$G$, and
non-trivial. Therefore, $G$ is not a twisted~product.

Let~$L$ be a divisible abelian group containing~$B$; to fix ideas, let
$L=\Q/\Z$, and let $j\colon B\to L$ be given by $j([x,y])={1\over p}$.

We first let~$f=j$ be defined only on~$B$, and we want to extend~$f$
to all of~$G$. First we extend~$f$ to~$x$. Since $x^p$ is the smallest
positive power of~$x$ that lies in~$B$, we choose a $p$-th root of
$$f(x^p)-{p\choose 2}j\left(\beta(x,x)\right)=0.$$
The easiest choice to make is $f(x)=0$, and we define~$f$
accordingly. Next we extend $f$ to~$y$; again, the least positive
power of~$y$ lying in $\langle B,x\rangle$ is~$y^p$, so we choose a
$p$-th root~of
$$f(y^p)-{p\choose 2}j\left(\beta(y,y)\right)=0,$$
so again we make the easiest choice and define $f(y)=0$.

Finally, we extend~$f$ to~$z$. Here we have that the least power
of~$z$ which lies in $\langle B,x,y\rangle$ is $z^p$, so we pick a
$p$-th root of
$$\eqalign{f(z^p)-{p\choose 2}j\left(\beta(z,z)\right)&= f(z^p)\cr
&=f([x,y])\cr
&= {1\over p}.\cr}$$
So we define $f(z)={1\over p^2}$, and extend to all of~$G$
using the formula in \ref{powerform}. 

The map $\phi\colon G\to G^{\hbox{\smallheadfont
ab}}\twist{j\circ\beta} L$ given
by $\phi(g)=(\overline{g},f(g))$ is the desired inclusion of~$G$ into
a twisted~product.

Note that we can embed~$G$ into a simpler twisted product by
replacing~$L$ with the image of~$f$, namely a cyclic group of
order~$p^2$.

{\bf Example \numbeq{zmodpsquare}.} For another example, consider
$\Z/p^2\Z$ as an extension of $\Z/p\Z$ by $\Z/p\Z$. We take $\Z/p^2\Z$
as the integers modulo~$p^2$, so~$G$ has a
distinguished generator, namely~$\overline{1}$, the class of~$1$. Thus
$B$ is the subgroup generated by $\overline{p}$. Here the map
\hbox{$(x,y)\mapsto [x,y]$} is trivial, so the universal map~$\beta$
associated to this bilinear map is symmetric. Again, let~$L=\Q/\Z$,
and define $f\colon B\to L$ by $f(\overline{p})={1\over p}$. 

To extend~$f$ to all of~$G$, we pick $\overline{1}$ as an element
of~$G$ not in~$B$; we need a $p$-th root of
$$f(\overline{p})-{p\choose 2}\beta(\overline{1},\overline{1})=
f(\overline{p}) = {1\over p}$$
so we let $f(\overline{1})={1\over p^2}$, and extend to all of~$G$
using the formula in \ref{powerform}. 

The map we obtain, $\varphi\colon \Z/p^2\Z \to \Z/p\Z\twist\beta\Q/\Z$ is
given by
$$\varphi(\overline{1}) = \left(\overline{1},{1\over p^2}\right).$$
Note that this embedds~$G$ into the second coordinate, and that the
first coordinated of the twisted product provides no new~information;
also note that in fact $$\Z/p\Z\twist\beta\Q/\Z\cong
\Z/p\Z\oplus\Z/\Z.$$ 

\rmrk{gettingridofdivsibility} In \ref{centextthentwist}, we use the
divisibility of~$L$ to obtain the map 
$\chi\colon C\to L$. It is also needed, implicitly, to extend~$f$
using \ref{fullextension}. It is possible to relax the divisibility
requirement given a specific group~$A$; but then we lose independence
from~$A$. For example, as noted in \ref{divisibility}, if~$A$ has
exponent~$k$ then we do not need the full divisibility of~$L$ to apply
\ref{fullextension}. To bypass the use of divisibility of~$L$ in
\ref{centextthentwist}, we proceed as follows: instead of using a
divisible group containing~$B$, we use a group containing the
universal object~$C$, and which has all the necessary roots to make
the argument work; thus, for example, if as above~$A$ has exponent
$k$, then we let~$L$ be any group containing~$C$ and having all~$k$-th
roots.

However, this process makes~$L$ depend not only on~$B$ and~$A$, but
also on~$G$ (since~$L$ depends now on~$C$, which in turn depends on
the commutator mapping, and hence on~$G$). Therefore, we do the
following: for each factor set $\gamma\in H^2(A,B)$, we let
$L_{\gamma}$ be the corresponding group obtained as in the paragraph
above. Finally, let~$L$ be the pushout of the $L_{\gamma}$ with
respect to the inclusions $B\hookrightarrow L_{\gamma}$, obtained by
embedding $B$ to~$C$ via $i_B$, and then~$C$ into~$L_{\gamma}$.  This
makes $L$ depend only on~$H^2(A,B)$, and therefore only on~$A$
and~$B$.

\Section{Cohomologically speaking}{cohomology}

The language of Cohomology of Groups is probably the best one to
express the result of \ref{centextandtwists}. It will also yield
further refinements and information about the given~construction.

\lemma{embedcentcohom}{Let $1\,\mapright{}\,B\,\mapright{}\,
G\,\mapright{}\, A\,\mapright{}\,1$ be a central extension of groups,
and let $\gamma$ in~$H^2(A,B)$ be a factor set representing this
extension, with associated transversal \hbox{$\ell\colon A\to G$}. Then there
exists an abelian group $L$ and an embedding
$j\colon B\to L$ such that $j^*(\gamma)\in H^2_{\rm
Bil}(A,L)$. Moreover, the group~$L$ and the embedding~$j$ can be
taken to depend only on~$B$.}

\proof Let $L$, $j$ and $f\colon G\to L$ be as constructed in
\ref{centextthentwist}. It suffices to show that $j\circ\gamma$ is
equal to the $\beta$ of the twisted product obtained in
\ref{centextthentwist}, as elements of~$H^2(A,L)$. Again we abuse
notation and use~$\beta$ to denote both the bilinear map from~$A$
to~$L$, and the induced bilinear map from~$G$ to~$L$. 

Let $h\colon A\to L$ be the set map given by $h(x)=f(\ell(x))$ for
every $x\in A$. Then $$h(e)=f(\ell(e))=f(e)=0,$$
and for any $x$ and~$y$ in $A$, we have
$$\eqalignno{\scriptstyle h(x)+h(y) - h(xy)&
\scriptstyle = f(\ell(x)) + f(\ell(y)) -
f(\ell(xy))\cr
&\scriptstyle =f(\ell(x)) + f(\ell(y)) - f(\ell(x)+ \ell(y)-
\gamma(x,y))\cr
&\scriptstyle =f(\ell(x)) + f(\ell(y)) - f(\ell(x) + \ell(y)) 
- f(-\gamma(x,y)) + \beta(\ell(x)+\ell (y),\gamma(x,y))\cr
&\scriptstyle =f(\ell(x)) + f(\ell(y)) -
f(\ell(x)+\ell(y)) + j(\gamma(x,y))\cr 
&\qquad\hbox{\smallheadfont(Since $\scriptstyle \gamma(x,y)$ lies in
$\scriptstyle B$, hence in the kernel
on the right of $\scriptstyle \beta$)}\cr
&\scriptstyle =f(\ell(x)) + f(\ell(y)) - f(\ell(x)) - f(\ell(y))
- \beta(\ell(x),\ell(y)) + j(\gamma(x,y))\cr
&\scriptstyle= j(\gamma(x,y)) - \beta(\ell(x),\ell(y))\cr
&\scriptstyle =j(\gamma(x,y))- \beta(x,y).\cr
&\qquad\hbox{\smallheadfont (Since $\scriptstyle \ell(x)$ is just the
image of~$\scriptstyle x$)}\cr}$$

So $\beta$ and $j\circ\gamma$ are equal in $H^2(A,L)$, as~claimed.\endproof

\thm{onlyBmatters}{Let $B$ be an abelian group. Then there
exists an abelian group~$L$, and an embedding $j\colon B\to L$
such that 
$$j^*(H^2(A,B))\subseteq H^2_{\rm Bil}(A,L)$$
for every abelian group~$A$.\noproof}

As noted before, $L$ will be a divisible abelian
group by construction, with the possibility of weakening this
requirement for a particular given~$A$.

Example \ref{zmodpsquare} calls our attention to the fact that the
embedding obtained may very well be ``uninteresting'', that is, that
it may be an embedding of our given group~$G$ into a direct sum (the
zero element of $H^2(A,L)$). So the obvious question to ask is, under
our construction, when does a
central extension get embedded into a direct sum?

Let $1\mapright{} B\mapright{j}L\mapright{p}L/B\mapright{} 1$
be the exact sequence associated to the inclusion of $B$ into
a divisible abelian group $L$. We will now use several results from
cohomology of groups. We refer the reader to either {\bf
[\cite{brown}]} or {\bf [\cite{weibel}]} for the details.

We have for every abelian group~$A$ a long exact
sequence of co\-ho\-mology~groups
$$\displaylines{0\mapright{}H^0(A,B)\mapright{j^*}H^0(A,L)\mapright{p^*}
H^0(A,L/B)\mapright{\delta_0}\hfill\cr
\hfill\mapright{\delta_0}H^1(A,B)\mapright{j^*}H^1(A,L)
\mapright{p^*}H^1(A,L/B)\mapright{\delta_1}H^2(A,B)
\mapright{j^*}H^2(A,L)\mapright{p^*}\cdots\cr}$$

Recall that in our case, with $A$ acting trivially on~$B$, $H^0(A,B)$ 
is just~$B$, and that $H^1(A,B)$ is the
group of all group morphisms from~$A$ to~$B$. So we see that
the above exact sequence yields the exact sequence
$$\scriptstyle0\mapright{}B\mapright{j}
L\mapright{p}L/B\mapright{\delta _0}
{\rm Hom}(A,B)\mapright{j^*}{\rm Hom}(A,L)\mapright{p*}{\rm
Hom}(A,L/B)
\mapright{\delta_1}H^2(A,B)\mapright{j^*}H^2(A,L)\mapright{p^*}\cdots$$
where ${\rm Hom}(A,B)$ is the set of group maps from~$A$ to~$B$,
etc. Since $p$ is surjective, $\delta_0$ is the zero map, and we have
the exact sequence
$$0\mapright{}{\rm Hom}(A,B)\mapright{j^*}{\rm
Hom}(A,L)\mapright{p*}{\rm Hom}(A,L/B)\mapright{\delta_1}
H^2(A,B)\mapright{j^*}H^2(A,L)\mapright{p^*}\cdots
$$

Also from Homological Algebra, we know that the right derived functor
of ${\rm Hom}$ is ${\rm Ext}$, that is, that we have an exact
sequence
$$0\mapright{}{\rm Hom}(A,B)\mapright{j^*}{\rm
Hom}(A,L)\mapright{p*}{\rm Hom}(A,L/B)\mapright{\delta_1}
{\rm Ext}(A,B)\mapright{j^*}{\rm Ext}(A,L)\mapright{p^*}\cdots$$

Since $L$ is divisible, it is an injective abelian group, so ${\rm
Ext}(A,L)=0$; hence $\delta_1$ is surjective onto ${\rm Ext}(A,B)$,
which as we know is a subgroup of $H^2(A,B)$. We conclude that the
kernel of $j^*\colon H^2(A,B) \to H^2(A,L),$
which is the map in \ref{onlyBmatters}, is none other than ${\rm
Ext}(A,B)$, that is, the abelian extensions of $B$ by $A$.

So we have proven:

\thm{kernelisext}{Let $B$ be an abelian group. Then there exists
an abelian group~$L$, and an embedding $j\colon B\to L$, such that
for any abelian group~$A$, 
$$j^*\colon H^2(A,B)\longrightarrow H^2(A,L)$$ 
has image contained in $H^2_{\rm Bil}(A,L)$, and kernel equal to ${\rm
Ext}(A,B)$. In other words, for every central extension $G$ of~$B$
by~$A$ there exists a twisted product $A\,\twist\beta\,L$ and an
embedding $G\hookrightarrow A\twist\beta L$, whose first component is
the projection of $G$ onto $A$. This twisted product is equivalent to
a direct sum of~$L$ and~$A$ if and only if $G$ is abelian.\noproof}

\cor{ifBisdiv}{Let $A$ be any abelian group, and
let~$B$ be divisible. Then every element of $H^2(A,B)$ is cohomologous
to a bilinear factor set. That is, every central extension of~$B$ by~$A$
is equivalent to a twisted product of~$B$ by~$A$.}

\proof We can let $L$ in \ref{centextthentwist} be equal to~$B$, and
$j\colon B\to L$ be the identity of~$B$. Then \ref{kernelisext} says
that ${\rm id}_B^*( = {\rm id}_{H^2(A,B)})$ has image in $H^2_{\rm
Bil}(A,B)$ and trivial kernel; so we have that $$H^2(A,B) = H^2_{\rm
Bil}(A,B)$$ as claimed.\endproof 

\thm{Afree}{If $A$ is a free abelian group, and~$B$ is any abelian
group, then $$H^2(A,B) = H^2_{\rm Bil}(A,B).$$}

\proof Let ${x_i}_{i\in I}$ be a basis for~$A$. It is not hard to
verify that we can apply the same process as in
Example~\ref{niltwoproduct} sequentially to the elements $x_i$. It is
not hard to verify that we will always be in Case~1 of the proof of
\ref{fullextension}, and thus that we may bypass the need for a
divisible group in that proof.\endproof

%
\ifnum0<\citations{\par\bigbreak
\filbreak{\bf References}\par\frenchspacing}\fi
%
\ifundefined{xthreeNB}\else
\item{\bf [\refer{threeNB}]}{Baumslag, G{.}, Neumann, B{.}H{.},
Neumann, H{.}, and Neumann, P{.}M. {\it On varieties generated by a
finitely generated group.\/} {\sl Math.\ Z.} {\bf 86} (1964)
pp.~\hbox{93--122}. {MR:30\#138}}\par\filbreak\fi
\ifundefined{xbergman}\else
\item{\bf [\refer{bergman}]}{Bergman, George M. {\it An Invitation to
General Algebra and Universal Constructions.\/} {\sl Berkeley
Mathematics Lecture Notes 7\/} (1995).}\par\filbreak\fi
\ifundefined{xordersberg}\else
\item{\bf [\refer{ordersberg}]}{Bergman, George M. {\it Ordering
coproducts of groups and semigroups.\/} {\sl J. Algebra 133\/} (1990)
no. 2, pp.~\hbox{313--339}. {MR:91j:06035}}\par\filbreak\fi
\ifundefined{xbrown}\else
\item{\bf [\refer{brown}]}{Brown, Kenneth S. {\it Cohomology of
Groups, 2nd Edition.\/} {\sl Graduate texts in mathematics 87\/},
Springer Verlag,~1994. {MR:96a:20072}}\par\filbreak\fi
\ifundefined{xmetab}\else
\item{\bf [\refer{metab}]}{Golovin, O. N. {\it Metabelian products of
groups.\/}
{\sl American Mathematical Society Translations}, series 2, {\bf 2} (1956),
pp.~\hbox{117--131.} {MR:17,824b}}\par\filbreak\fi
\ifundefined{xhall}\else
\item{\bf [\refer{hall}]}{Hall, M. {\it The Theory of Groups.\/}
Mac~Millan Company,~1959. {MR:21\#1996}}\par\filbreak\fi
\ifundefined{xphall}\else
\item{\bf [\refer{phall}]}{Hall, P. {\it Verbal and marginal
subgroups.} {\sl J.\ Reine\ Angew.\ Math.\/} {\bf 182} (1940)
pp.~\hbox{156--157.} {MR:2,125i}}\par\filbreak\fi
\ifundefined{xheineken}\else
\item{\bf [\refer{heineken}]}{Heineken, H. {\it Engelsche Elemente der
L\"ange drei,\/} {\sl Illinois Journal of Math.} {\bf 5} (1961)
pp.~\hbox{681--707.} {MR:24\#A1319}}\par\filbreak\fi
\ifundefined{xherman}\else
\item{\bf [\refer{herman}]}{Herman, Krzysztof. {\it Some remarks on
the twelfth problem of Hanna Neumann.\/} {\sl Publ.\ Math.\ Debrecen}
{\bf 37} (1990)  no. 1--2, pp.~\hbox{25--31.} {MR:91f:20030}}\par\filbreak\fi
\ifundefined{xepisandamalgs}\else
\item{\bf [\refer{episandamalgs}]}{Higgins, Peter M. {\it Epimorphisms
and amalgams.} {\sl
Colloq.\ Math.} {\bf 56} no.~1 (1988) pp.~\hbox{1--17.}
{MR:89m:20083}}\par\filbreak\fi
\ifundefined{xhigmanpgroups}\else
\item{\bf [\refer{higmanpgroups}]}{Higman, Graham. {\it Amalgams of
$p$-groups.\/} {\sl J. of~Algebra} {\bf 1} (1964)
pp.~\hbox{301--305.} {MR:29\#4799}}\par\filbreak\fi
\ifundefined{xhigmanremarks}\else
\item{\bf [\refer{higmanremarks}]}{Higman, Graham. {\it Some remarks
on varieties of groups.\/} {\sl Quart.\ J.\ of Math.\ (Oxford) (2)} {\bf
10} (1959), pp.~\hbox{165--178.} {MR:22\#4756}}\par\filbreak\fi
\ifundefined{xhughes}\else
\item{\bf [\refer{hughes}]}{Hughes, N.J.S. {\it The use of bilinear
mappings in the classification of groups of class~$2$.\/} {\sl Proc.\
Amer.\ Math.\ Soc.\ } {\bf 2} (1951) pp.~\hbox{742--747.}
{MR:13,528e}}\par\filbreak\fi
\ifundefined{xisbelltwo}\else
\item{\bf [\refer{isbelltwo}]}{Howie, J.~M., Isbell, J.~R. {\it
Epimorphisms and dominions II.\/} {\sl Journal of Algebra {\bf
6}}(1967) pp.~\hbox{7--21.} {MR:35\#105b}}\par\filbreak\fi
\ifundefined{xisaacs}\else
\item{\bf [\refer{isaacs}]}{Isaacs, I.M., Navarro, Gabriel. {\it
Coprime actions, fixed-point subgroups and irreducible induced
characters.} {\sl J.~of Algebra} {\bf 185} (1996) no.~1,
pp.~\hbox{125--143.} {MR:97g:20009}}\par\filbreak\fi
\ifundefined{xisbellone}\else
\item{\bf [\refer{isbellone}]}{Isbell, J. R. {\it Epimorphisms and
dominions} in {\sl 
Proc.~of the Conference on Categorical Algebra, La Jolla 1965,\/}
pp.~\hbox{232--246.} Lange and Springer, New
York~1966. {MR:35\#105a}(The statement of the
Zigzag Lemma for {\it rings} in this paper is incorrect. The correct
version is stated in~{\bf [\cite{isbellfour}]})}\par\filbreak\fi
\ifundefined{xisbellthree}\else
\item{\bf [\refer{isbellthree}]}{Isbell, J. R. {\it Epimorphisms and
dominions III.} {\sl Amer.\ J.\ Math.\ }{\bf 90} (1968)
pp.~\hbox{1025--1030.} {MR:38\#5877}}\par\filbreak\fi
\ifundefined{xisbellfour}\else
\item{\bf [\refer{isbellfour}]}{Isbell, J. R. {\it Epimorphisms and
dominions IV.} {\sl Journal\ London Math.\ Society~(2),}
{\bf 1} (1969) pp.~\hbox{265--273.} {MR:41\#1774}}\par\filbreak\fi
\ifundefined{xjones}\else
\item{\bf [\refer{jones}]}{Jones, Gareth A. {\it Varieties and simple
groups.\/} {\sl J.\ Austral.\ Math.\ Soc.\ } {\bf 17}
pp.~\hbox{163--173.} {MR:49\#9081}}\par\filbreak\fi
\ifundefined{xjonsson}\else
\item{\bf [\refer{jonsson}]}{J\'onsson, B. {\it Varieties of groups of
nilpotency three.} {\sl Notices Amer.\ Math.\ Soc.} {\bf 13} (1966)
pp.~488.}\par\filbreak\fi
\ifundefined{xwreathext}\else
\item{\bf [\refer{wreathext}]}{Kaloujnine, L. and Krasner, Marc. {\it
Produit complet des groupes de permutations et le probl\`eme
d'extension des groupes III.} {\sl Acta Sci.\ Math.\ Szeged} {\bf 14}
(1951) pp.~\hbox{69--82}. {MR:14,242d}}\par\filbreak\fi
\ifundefined{xkhukhro}\else
\item{\bf [\refer{khukhro}]}{Khukhro, Evgenii I. {\it Nilpotent Groups
and their Automorphisms.} {\sl de Gruyter Expositions in Mathematics}
{\bf 8}, New York 1993. {MR:94g:20046}}\par\filbreak\fi
\ifundefined{xkleimanbig}\else
\item{\bf [\refer{kleimanbig}]}{Kle\u{\i}man, Yu.~G. {\it On
identities in groups.\/} {\sl Trans.\ Moscow Math.\ Soc.\ } 1983,
Issue 2, pp.~\hbox{63--110}. {MR:84e:20040}}\par\filbreak\fi
\ifundefined{xthirtynine}\else
\item{\bf [\refer{thirtynine}]}{Kov\'acs, L.~G. {\it The thirty-nine
varieties.} {\sl Math.\ Scientist} {\bf 4} (1979)
pp.~\hbox{113--128.} {MR:81m:20037}}\par\filbreak\fi
\ifundefined{xlevione}\else
\item{\bf [\refer{levione}]}{Levi, F.~W. {\it Groups on which the
commutator relation 
satisfies certain algebraic conditions.\/} {\sl J.\ Indian Math.\ Soc.\ New
Series} {\bf 6}(1942), pp.~\hbox{87--97.} {MR:4,133i}}\par\filbreak\fi
\ifundefined{xgermanlevi}\else
\item{\bf [\refer{germanlevi}]}{Levi, F.~W. and van der Waerden,
B.~L. {\it \"Uber eine 
besondere Klasse von Gruppen.\/} {\sl Abhandl.\ Math.\ Sem.\ Univ.\ Hamburg}
{\bf 9}(1932), pp.~\hbox{154--158.}}\par\filbreak\fi
\ifundefined{xlichtman}\else
\item{\bf [\refer{lichtman}]}{Lichtman, A.~L. {\it Necessary and
sufficient conditions for the residual nilpotence of free products of
groups.\/} {\sl J. Pure and Applied Algebra} {\bf 12} no. 1 (1978),
pp.~\hbox{49--64.} {MR:58\#5938}}\par\filbreak\fi
\ifundefined{xmaxofan}\else
\item{\bf [\refer{maxofan}]}{Liebeck, Martin W.; Praeger, Cheryl E.;
and Saxl, Jan. {\it A classification of the maximal subgroups of the
finite alternating and symmetric groups.\/} {\sl J. of Algebra} {\bf
111}(1987), pp.~\hbox{365--383.} {MR:89b:20008}}\par\filbreak\fi
\ifundefined{xepisingroups}\else
\item{\bf [\refer{episingroups}]}{Linderholm, C.E. {\it A group
epimorphism is surjective.\/} {\sl Amer.\ Math.\ Monthly\ }77
pp.~\hbox{176--177.}}\par\filbreak\fi
\ifundefined{xmckay}\else
\item{\bf [\refer{mckay}]}{McKay, Susan. {\it Surjective epimorphisms
in classes
of groups.} {\sl Quart.\ J.\ Math.\ Oxford (2),\/} {\bf 20} (1969),
pp.~\hbox{87--90.} {MR:39\#1558}}\par\filbreak\fi
\ifundefined{xmaclane}\else
\item{\bf [\refer{maclane}]}{Mac Lane, Saunders. {\it Categories for
the Working Mathematician.} {\sl Graduate texts in mathematics 5},
Springer Verlag (1971). {MR:50\#7275}}\par\filbreak\fi
\ifundefined{xbilinear}\else
\item{\bf [\refer{bilinear}]}{Magidin, A. {\it Bilinear maps and 
2-nilpotent groups.\/} August 1996, 7~pp.}\par\filbreak\fi
\ifundefined{xbilinearprelim}\else
\item{\bf [\refer{bilinearprelim}]}{Magidin, Arturo. {\it Bilinear maps
and central extensions of abelian groups.\/} In~preparation.}\par\filbreak\fi
\ifundefined{xprodvar}\else
\item{\bf [\refer{prodvar}]}{Magidin, A. {\it Dominions in product
varieties of groups.\/} May 1997, 21~pp.}\par\filbreak\fi
\ifundefined{xnildoms}\else
\item{\bf [\refer{nildoms}]}{Magidin, A. {\it Dominions in varieties
of nilpotent groups.\/} December 1996, 27~pp.}\par\filbreak\fi
\ifundefined{xnildomsprelim}\else
\item{\bf [\refer{nildomsprelim}]}{Magidin, Arturo. {\it Dominions in
varieties of nilpotent groups.\/} In~preparation.}\par\filbreak\fi
\ifundefined{xntwodoms}\else
\item{\bf [\refer{ntwodoms}]}{Magidin, A. {\it Dominions in the variety of
2-nilpotent groups.\/} May 1996, 6~pp.}\par\filbreak\fi
\ifundefined{xfgnilgroups}\else
\item{\bf [\refer{fgnilgroups}]}{Magidin, A. {\it Dominions of
finitely generated nilpotent groups.\/} October~1997,
10~pp.}\par\filbreak\fi
\ifundefined{xepis}\else
\item{\bf [\refer{epis}]}{Magidin, A. {\it Non-surjective epimorphisms
in varieties of groups and other results.\/} February 1997,
13~pp.}\par\filbreak\fi
\ifundefined{xoddsandends}\else
\item{\bf [\refer{oddsandends}]}{Magidin, A. {\it Some odds and
ends.\/} June 1996, 3~pp.}\par\filbreak\fi
\ifundefined{xpropdom}\else
\item{\bf [\refer{propdom}]}{Magidin, A. {\it Some properties of
dominions in varieties of groups.\/} March 1997, 13~pp.}\par\filbreak\fi
\ifundefined{xmagnus}\else
\item{\bf [\refer{magnus}]}{Magnus, Wilhelm; Karras, Abraham; and
Solitar, Donald. {\it Combinatorial Group Theory.\/} 2nd Edition; Dover
Publications, Inc.~1976. {MR:53\#10423}}\par\filbreak\fi
\ifundefined{xamalgtwo}\else
\item{\bf [\refer{amalgtwo}]}{Maier, Berthold J. {\it Amalgame
nilpotenter Gruppen
der Klasse zwei II.\/} {\sl Publ.\ Math.\ Debrecen} {\bf 33}(1986),
pp.~\hbox{43--52.} {MR:87k:20050}}\par\filbreak\fi
\ifundefined{xnilexpp}\else
\item{\bf [\refer{nilexpp}]}{Maier, Berthold J. {\it On nilpotent
groups of exponent $p$.\/} {\sl Journal of~Algebra} {\bf 127} (1989)
pp.~\hbox{279--289.} {MR:91b:20046}}\par\filbreak\fi
\ifundefined{xmaltsev}\else
\item{\bf [\refer{maltsev}]}{Maltsev, A.~I. {\it Generalized
nilpotent algebras and their associated groups.} (Russian) {\sl
Mat.\ Sbornik N.S.} {\bf 25(67)} (1949) pp.~\hbox{347--366.} ({\sl
Amer.\ Math.\ Soc.\ Translations Series 2} {\bf 69} 1968,
pp.~\hbox{1--21.}) {MR:11,323b}}\par\filbreak\fi
\ifundefined{xmaltsevtwo}\else
\item{\bf [\refer{maltsevtwo}]}{Maltsev, A.~I. {\it Homomorphisms onto
finite groups.} (Russian) {\sl Ivanov. gosudarst. ped. Inst., u\v
cenye zap., fiz-mat. Nuak} {\bf 18} (1958)
\hbox{pp. 49--60.}}\par\filbreak\fi
\ifundefined{xmorandual}\else
\item{\bf [\refer{morandual}]}{Moran, S. {\it Duals of a verbal
subgroup.\/} {\sl J.\ London Math.\ Soc.} {\bf 33} (1958)
pp.~\hbox{220--236.} {MR:20\#3909}}\par\filbreak\fi
\ifundefined{xhneumann}\else
\item{\bf [\refer{hneumann}]}{Neumann, Hanna. {\it Varieties of
Groups.\/} {\sl Ergebnisse der Mathematik und ihrer Grenz\-ge\-biete\/}
New series, Vol.~37, Springer Verlag~1967. {MR:35\#6734}}\par\filbreak\fi
\ifundefined{xneumannwreath}\else
\item{\bf [\refer{neumannwreath}]}{Neumann, Peter M. {\it On the
structure of standard wreath products of groups.\/} {\sl Math.\
Zeitschr.\ }{\bf 84} (1964) pp.~\hbox{343--373.} {MR:32#5719}}\par\filbreak\fi
\ifundefined{xpneumann}\else
\item{\bf [\refer{pneumann}]}{Neumann, Peter M. {\it Splitting groups
and projectives
in varieties of groups.\/} {\sl Quart.\ J.\ Math.\ Oxford} (2), {\bf
18} (1967),
pp.~\hbox{325--332.} {MR:36\#3859}}\par\filbreak\fi
\ifundefined{xoates}\else
\item{\bf [\refer{oates}]}{Oates, Sheila. {\it Identical Relations in
Groups.\/} {\sl J.\ London Math.\ Soc.} {\bf 38} (1963),
pp.~\hbox{71--78.} {MR:26\#5043}}\par\filbreak\fi
\ifundefined{xolsanskii}\else
\item{\bf [\refer{olsanskii}]}{Ol'\v{s}anski\v{\i}, A. Ju. {\it On the
problem of a finite basis of identities in groups.\/} {\sl
Izv.\ Akad.\ Nauk.\ SSSR} {\bf 4} (1970) no. 2
pp.~\hbox{381--389.}}\par\filbreak\fi
\ifundefined{xremak}\else
\item{\bf [\refer{remak}]}{Remak, R. {\it \"Uber minimale invariante
Untergruppen in der Theorie der endlichen Gruppen.\/} {\sl
J.\ reine.\ angew.\ Math.} {\bf 162} (1930),
pp.~\hbox{1--16.}}\par\filbreak\fi
\ifundefined{xclassifthree}\else
\item{\bf [\refer{classifthree}]}{Remeslennikov, V. N. {\it Two
remarks on 3-step nilpotent groups} (Russian) {\sl Algebra i Logika
Sem.} (1965) no.~2 pp.~\hbox{59--65.} {MR:31\#4838}}\par\filbreak\fi
\ifundefined{xrotman}\else
\item{\bf [\refer{rotman}]}{Rotman, J.J. {\it Introduction to the Theory of
Groups}, 4th edition. {\sl Graduate texts in mathematics 119},
Springer Verlag,~1994. {MR:95m:20001}}\par\filbreak\fi
\ifundefined{xsaracino}\else
\item{\bf [\refer{saracino}]}{Saracino, D. {\it Amalgamation bases for
nil-$2$ groups.\/} {\sl Alg.\ Universalis\/} {\bf 16} (1983),
pp.~\hbox{47--62.} {MR:84i:20035}}\par\filbreak\fi
\ifundefined{xscott}\else
\item{\bf [\refer{scott}]}{Scott, W.R. {\it Group Theory.} Prentice
Hall,~1964. {MR:29\#4785}}\par\filbreak\fi
\ifundefined{xsmelkin}\else
\item{\bf [\refer{smelkin}]}{Smel'kin, A.L. {\it Wreath products and
varieties of groups} [Russian] {\sl Dokl.\ Akad.\ Nauk S.S.S.R.\/} {\bf
157} (1964), pp.~\hbox{1063--1065} Transl.: {\sl Soviet Math.\ Dokl.\ } {\bf
5} (1964), pp.~\hbox{1099--1011}. {MR:33\#1352}}\par\filbreak\fi
\ifundefined{xstruiktwo}\else
\item{\bf [\refer{struiktwo}]}{Struik, Ruth Rebekka. {\it On nilpotent
products of cyclic groups II.\/} {\sl Canadian Journal of
Mathematics\/} {\bf 13} (1961) pp.~\hbox{557--568.}
{MR:26\#2486}}\par\filbreak\fi
\ifundefined{xvlee}\else
\item{\bf [\refer{vlee}]}{Vaughan-Lee, M{.} R{.} {\it Uncountably many
varieties of groups.\/} {\sl Bull.\ London Math.\ Soc.} {\bf 2} (1970)
pp.~\hbox{280--286.} {MR:43\#2054}}\par\filbreak\fi
\ifundefined{xweibel}\else
\item{\bf [\refer{weibel}]}{Weibel, Charles. {\it Introduction to
Homological Algebra.\/} Cambridge University
Press~1994. {MR:95f:18001}}\par\filbreak\fi 
\ifundefined{xweigelone}\else
\item{\bf [\refer{weigelone}]}{Weigel, T.S. {\it Residual properties
of free groups.\/} {\sl J.\ of Algebra} {\bf 160} (1993)
pp.~\hbox{14--41.} {MR:94f:20058a}}\par\filbreak\fi
\ifundefined{xweigeltwo}\else
\item{\bf [\refer{weigeltwo}]}{Weigel, T.S. {\it Residual properties
of free groups II.\/} {\sl Comm.\ in Algebra} {\bf 20}(5) (1992)
pp.~\hbox{1395--1425.} {MR:94f:20058b}}\par\filbreak\fi
\ifundefined{xweigelthree}\else 
\item{\bf [\refer{weigelthree}]}{Weigel, T.S. {\it Residual Properties
of free groups III.\/} {\sl Israel J.\ Math.\ } {\bf 77} (1992)
pp.~\hbox{65--81.} {MR:94f:20058c}}\par\filbreak\fi
\ifundefined{xzstwo}\else
\item{\bf [\refer{zstwo}]}{Zariski, Oscar and Samuel, Pierre. {\it
Commutative Algebra}, Volume
II. Springer-Verlag~1976. {MR:52\#10706}}\par\filbreak\fi
\ifnum0<\citations\nonfrenchspacing\fi
\smallskip
{\frenchspacing
\obeylines\baselineskip=12pt plus3pt minus1pt
\it
Arturo Magidin
Department of Mathematics
970 Evans Hall
University of California at Berkeley
Berkeley, CA 94720, USA
e-mail: magidin@math.berkeley.edu\par
}

\vfill
\eject
\immediate\closeout\aux
\end